\documentclass[preprint,12pt]{elsarticle}
\usepackage{amssymb}
\usepackage[para,online,flushleft]{threeparttable}
\graphicspath{ {./figures/} }
\usepackage{hyperref}
\usepackage{float}
\usepackage{verbatim}
\restylefloat{figure}
\restylefloat{table}

\journal{Iran Journal of Computer Science}

\usepackage{titlesec}
\titleclass{\subsubsubsection}{straight}[\subsection]
\newcounter{subsubsubsection}[subsubsection]
\renewcommand\thesubsubsubsection{\thesubsubsection.\arabic{subsubsubsection}}
\titleformat{\subsubsubsection}{\normalfont\normalsize\itshape}{\thesubsubsubsection.\space}{0em}{}
\titlespacing*{\subsubsubsection}{0pt}{2ex plus 1ex minus .2ex}{0.75ex plus .2ex}
\makeatletter
\def\toclevel@subsubsubsection{4}
\def\l@subsubsubsection{\@dottedtocline{4}{7em}{4em}}
\makeatother

\usepackage[table]{xcolor}
\usepackage{graphicx}
\usepackage{amsmath,amssymb,amsfonts}
\usepackage{algorithmic}
\usepackage{graphicx}
\usepackage{textcomp}
\usepackage{booktabs}
\usepackage{xcolor}
\usepackage{array}
\usepackage{multirow}
\usepackage{url}
\usepackage{float}
\usepackage{pifont}

\floatstyle{plaintop}
\restylefloat{table}
\usepackage{caption}
\usepackage{subcaption}

\newcommand{\blt}[0]{\color{black}}

\usepackage{booktabs}

\newif\ifblackandwhite
\blackandwhitetrue
\usepackage{pdflscape}
\usepackage{colortbl}

\usepackage{booktabs}

\usepackage{adjustbox}
\usepackage{rotating}

\def\BibTeX{{\rm B\kern-.05em{\sc i\kern-.025em b}\kern-.08em
    T\kern-.1667em\lower.7ex\hbox{E}\kern-.125emX}}

\usepackage{geometry}
\begin{document}
\begin{frontmatter}

\begin{titlepage}
\begin{center}
\vspace*{0.5cm}

\textbf{Integrating Random Regret Minimization-Based Discrete Choice Models with Mixed Integer Linear Programming for Revenue Optimization}
\vspace{2cm}

Amirreza Talebi$^{a,c}$ (talebi.14@osu.edu)\\
Sayed Pedram Haeri Boroujeni$^b$ (shaerib@g.clemson.edu)\\
Abolfazl Razi$^b$ (arazi@clemson.edu)\\

\hspace{10pt}

\begin{flushleft}
\small  
$^a$Department of Integrated Systems Engineering, The Ohio State University, Columbus, OH, USA \\[1mm]
$^b$School of Computing, Clemson University, Clemson, SC 29632, USA\\[1mm]
$^c$ENSIMAG, Grenoble Alpes University, Grenoble, France\\[1mm]

\vspace{2.5cm}

\textbf{Corresponding Author:} \\
Amirreza Talebi\\
Department of Integrated Systems Engineering, The Ohio State University, Columbus, OH, USA \\
Email: talebi.14@osu.edu\\

\end{flushleft}        
\end{center}
\end{titlepage}

\title{Integrating Random Regret Minimization-Based Discrete Choice Models with Mixed Integer Linear Programming for Revenue Optimization}


\author{Amirreza Talebi$^{a,c}$ Sayed Pedram Haeri Boroujeni$^{b}$, Abolfazl Razi$^{b}$}

\affiliation{organization={Department of Integrated Systems Engineering},
            addressline={The Ohio State University}, 
            city={Columbus},
            postcode={43210}, 
            state={OH},
            country={USA}}

\affiliation{organization={School of Computing},
            addressline={Clemson University}, 
            city={Clemson},
            postcode={29632}, 
            state={SC},
            country={USA}}

\affiliation{organization={ENSIMAG},
            addressline={Grenoble Alpes University}, 
            city={Grenoble},
            postcode={38400}, 
            state={Saint-Martin d'Hères},
            country={France}}

\begin{abstract}
This paper explores the critical domain of Revenue Management (RM) within Operations Research (OR), focusing on intricate pricing dynamics. Utilizing Mixed Integer Linear Programming (MILP) models, the study enhances revenue optimization by considering product prices as decision variables and emphasizing the interplay between demand and supply. Recent advancements in Discrete Choice Models (DCMs), particularly those rooted in Random Regret Minimization (RRM) theory, are investigated as potent alternatives to established Random Utility Maximization (RUM) based DCMs.
Despite the widespread use of DCMs in RM, a significant gap exists between cutting-edge RRM-based models and their practical integration into RM strategies. The study addresses this gap by incorporating an advanced RRM-based DCM into MILP models, addressing pricing challenges in both capacitated and uncapacitated supply scenarios. The developed models demonstrate the feasibility and offer diverse interpretations of consumer choice behavior, drawing inspiration from established RUM-based frameworks. This research contributes to bridging the existing gap in the application of advanced DCMs within practical RM implementations.
\end{abstract}


\end{frontmatter}

\allowdisplaybreaks

\section{Introduction}
Advancements in Operations Research (OR) Modeling, coupled with innovative solving methods, have efficiently addressed large-scale challenges across economics, transportation, supply chain, revenue management, etc. The primary goal of these models is to optimize revenue, minimize costs, or mitigate pollution. Consequently, sophisticated demand models are under development to quantify these factors precisely. Recognizing the intricate interplay between demand and supply, our approach integrates both aspects, emphasizing simultaneous consideration of the supply side during supply chain planning. Notably, the pricing decision involves suppliers setting prices while buyers determine their willingness to pay, exemplified by the varying affordability of a \$10 rose in different neighborhoods. This underscores the crucial need for sellers to comprehend and anticipate customer purchasing behavior. Concurrently, Revenue Management (RM), within the realm of OR, concentrates on maximizing revenue through strategic pricing and capacity allocation. Our work uniquely combines economic and OR perspectives, incorporating customer behavior modeling for optimized supply systems.

Historically, demand and revenue-maximizing prices were determined using simple linear demand functions. Additionally, the main assumption of consumer demand was that a customer would buy a product if it were available and would not buy at all if it were not available \cite{strauss2018review}. Moreover, the supply and demand aspects were treated separately, meaning the supplier would decide on the price and capacity of the products regardless of the demand, which could be inefficient as the demand and supply sides are intricately intertwined \cite{pacheco2016new}.

An essential assumption in RM is the variability of product prices for each customer. Precise understanding of customer behavior is necessary to optimize revenue. Discrete Choice Models (DCMs), which elucidate and predict customer choices, significantly contribute to customer behavior modeling. These models consider customer preferences, leading to efficient resource utilization and substantial revenue improvements. Numerous empirical studies underscore the impact of discrete choice modeling in diverse sectors such as RM \cite{talluri2004revenue}, labor market \cite{kornstad2007discrete}, food industry \cite{gracia2008demand}, and airline industry \cite{vulcano2010om}.

 Random Utility Maximization Theory (RUM), established by \cite{mcfadden1978modelling, mcfadden2000mixed}, forms the foundation for discrete choice modeling. The core assumption of RUM is that decision-makers are perfectly rational and aim to maximize their utilities. Recent studies, such as those by \cite{chorus2008random, chorus2010new, chorus2012rrm} propose an alternative approach by investigating DCMs derived from Random Regret Minimization (RRM) where the customer tries to minimize the regret caused by not choosing an alternative. Conditions exist where DCMs derived from either theory outperform the other, revealing diverse customer preferences \cite{hensher2013random, greene2011ru, li2018hybrid}.

While existing literature often examines demand and supply separately or investigates supply based on a simplistic demand formulation, a comprehensive approach is needed for effective revenue generation \cite{pacheco2016new}. Neglecting the interplay between supply and demand or focusing solely on one side can lead to suboptimal pricing strategies, resulting in reduced gains, the bullwhip effect, and discordance in supply chains. Recognizing that customer decisions drive demand underscores the vital importance of optimizing and coordinating both supply and demand concurrently. In this regard, \cite{pacheco2017integrating} innovatively crafted a framework that characterizes how demand responds to suppliers' decisions, linking supply and demand through price as a crucial variable. They linearized DCMs into predictors of demand using a Mixed Integer Linear Programming (MILP) model, emphasizing price variables in the objective function.

In related literature, other works have incorporated DCMs in OR to forecast demand and maximize revenue. For example, \cite{pacheco2016new, paneque2017integrating} addressed challenges arising from the non-linearity of DCM models in optimization. To address this, utility functions (RUM-based) were integrated into MILP constraints, introducing individual utility functions as scores to be maximized in discrete optimization, which is crucial for characterizing the demand side of the OR problem \cite{haeri2023hybrid, boroujeni2021novel}.

To date, all demand estimation models in the realm of Revenue Management (RM) have relied on RUM-based DCMs to maximize revenue. To the best of our knowledge, there are no MILP models based on RRM theory to maximize revenue. In this paper, we introduce a novel MILP model developed based on RRM theory. The main distinction between our proposed algorithm and others lies in the underlying assumption about customer behavior. In MILP models based on RUM theory, it is assumed that customers are utility maximizers. These models determine prices based on the expectation that customers will maximize their utilities while suppliers focus on revenue maximization.
Contrastingly, our algorithm operates on the premise that customers are regret minimizers. They select products aiming to minimize their regrets, while suppliers continue to strive for revenue maximization. This fundamental difference introduces novel challenges.
In RUM-based MILP models, customers compute the expected utility gained from purchasing each product and select the one that maximizes their utility. Conversely, in our model based on RRM, we compare all products pairwise. This approach enables customers to choose the product with the least regret. To do so, we have designed the MILP model with additional constraints.

Existing literature in RM predominantly applies RUM-based MILP to maximize revenue and address the pricing issue. However, in this paper, an RRM-based MILP model has been proposed for the first time. The main contributions of this work are as follows: 
\begin{itemize}
    \item In many real-world scenarios, customers prioritize minimizing regret when making purchasing decisions, rendering RUM-based MILP models inadequate for accurate pricing and revenue maximization. To address this limitation, we have developed a MILP model based on an advanced RRM-based DCM. This model is capable of incorporating customers' individual heterogeneity in revenue optimization and pricing strategies.
     \item Our proposed RRM-based model accommodates scenarios where suppliers have limited capacity to offer products to clients. Additionally, we have developed an alternative RRM-based model designed for scenarios where capacity is unlimited. 
    There are two main challenges in developing such models: \begin{itemize}
        \item The DCMs become highly non-linear and non-convex when logit functions are used to estimate the probability of a customer choosing a product, as demonstrated by \cite{pacheco2016new}. To overcome this challenge, we developed a linear formulation based on RRM theory, directly capturing how a customer selects a product to minimize regret. Additionally, we incorporated Gumbel-distributed noise to accommodate the uncertainties in consumer choices, which is consistent with RRM theory. Moreover, within the framework of our MILP model, we determine the average maximum revenue by sampling from the Gumbel distribution multiple times, a common practice in stochastic optimization approaches. This approach mirrors the technique used for the RUM-based MILP developed by \cite{pacheco2016new}.
        \item The computation of regret, even without considering the non-linearity introduced by DCMs, as mentioned in the previous point, inherently involves both minimization and maximization when comparing products. To handle this challenge, we employed linearization techniques in our MILP models for both capacitated and uncapacitated scenarios. 
    \end{itemize}
    \item We conducted experiments to assess the computational performance and solution quality of our proposed models compared to their RUM-based counterparts developed by  \cite{pacheco2016new, pacheco2017integrating}. 
\end{itemize}
Finally, we demonstrate that when customers are regret minimizers, the maximum revenue achieved by our models differs significantly from that of RUM-based counterparts. This discrepancy can have a profound impact on the production plans of suppliers targeting regret-minimizing customers with their products. In particular, our results demonstrate that RUM-based MILP models can significantly cause a supplier to lose revenue when customers are regret minimizers.

The primary objective of this paper is to develop a robust framework for RM capable of maximizing revenue and optimizing product pricing for a supplier in scenarios with both capacitated and uncapacitated supply while accounting for customer choice behavior. Our hypothesis posits that customers are regret minimizers in certain scenarios, aiming to select products or opt out to minimize their regrets through pairwise comparisons of products. We contend that current frameworks relying on RUM theory may underestimate prices and lead to revenue losses.
To address this, we consider both capacitated and uncapacitated supply quantities within our framework. Then, we introduce the RRM-based behavior of customers' decision-making together with capacity, pricing, and revenue concerns to the MILP. Through experiments, we demonstrate that our framework based on RRM theory yields higher revenue by aligning with the behavior of customers under RRM theory.

For the dataset, we selected a range of numerical values representing potential product prices. Additionally, we incorporated randomly drawn samples from a Gumbel(0,1) distribution to represent the unobserved regrets and utilities of customers. To ensure robustness in our results, We generated multiple draws of such data, taking into account the stochastic nature of unobserved regrets and utilities, a methodology also employed by \cite{pacheco2016new, pacheco2017integrating}.

Our results initially demonstrate the effectiveness of our developed models for both capacitated and uncapacitated versions. It is evident that the constrained version involves a larger number of variables and constraints compared to its unconstrained counterpart, resulting in longer solving times. Furthermore, our model solutions indicate that the RUM-based approach leads to revenue loss when customers are regret minimizers. This underscores the significance of employing RRM-based MILPs for such scenarios, which results in higher revenue for suppliers while still meeting customer needs. However, our models inherently entail more variables and constraints compared to their RUM-based counterparts, given the intrinsic nature of RRM theory. This characteristic inevitably leads to longer solution times.

The paper is structured as follows:

In Section \ref{sec2}, a comprehensive literature review is provided, encompassing revenue management, microeconomic behavioral theories, discrete choice modeling, choice-based optimization framework, and data collection.
Section \ref{sec3} introduces a novel representation of demand and supply based on RRM through MILP formulation.
Section \ref{sec4} conducts experiments, and Section \ref{sec5} concludes. 

\section{Related Work}
\label{sec2}

Within RM, pricing decisions hold prominent importance among various aspects to increase revenue. Several RM techniques aim to leverage customer heterogeneity for a product set at a specific time. They often focus on fixing alternatives or products and time dimensions while exploring customer differences. In our context, we emphasize the significance of customer heterogeneity in purchasing behavior. If all customers perceive a product uniformly with identical purchasing behavior, the customer dimension diminishes, leading to unexploited variations in willingness to pay. These variations, influenced by socio-economic factors like marital status, salary, age, gender, etc., play a critical role in maximizing revenue. On the other hand, optimizing revenue hinges on accurate demand prediction, which is crucial for effective supply planning and pricing strategies. 
The demand function is influenced by various variables, encompassing customer heterogeneity, prices, etc \cite{talluri2006theory, SADRALASHRAFI201882}. Specifically, researchers have extensively focused on DCMs, particularly due to their capacity to predict demand by accounting for consumer heterogeneity in tastes and attributes \cite{talluri2006theory, mcfadden2000mixed, chorus2008random}. RUM-based DCMs assume decision-makers are rational utility maximizers and models derived from these theories aim to reveal influential factors shaping customer choice preferences (e.g., see \cite{kemperman2021review} for a recent review). However, in some cases, customers prefer regret minimization over utility maximization when making choices. Therefore, deriving DCMs from suitable behavioral theories is crucial \cite{hensher2013random, gusarov2020exploration}. In this regard, \cite{chorus2014random} introduced a regret-based discrete choice model for consumer preference analysis with multinomial choice sets. Their empirical assessment, based on 43 comparisons between RUM and RRM models in the literature, revealed the superiority of RRM-based DCMs over RUM-based ones. The study demonstrated that both RRM-based and hybrid RRM-RUM-based DCMs consistently outperformed their RUM-based counterparts in terms of prediction accuracy and overall model fit. The authors emphasized that although performance differences were minor, the models led to significant managerial implications and diverse interpretations of behavior. 

The reference \cite{masiero2019understanding} conducted a choice experiment to investigate tourists' hotel preferences. They created hypothetical hotel alternatives and randomized factors, including price, hotel rating, and walking time to key attractions, to form choice sets. RRM and RUM-based DCMs were utilized to estimate the coefficients of these factors, and they figured out the RRM-based DCMs performed significantly better than their RUM-based counterparts. Additionally, they showed that quantitative simulations can be performed to study the demand and pricing issues while incorporating customers’ heterogeneity into the model. 

The reference \cite{sharma2019park} investigates the factors influencing park-and-ride (PNR) lot choice behavior among users. Two models, based on RUM and RRM concepts, are developed to understand why users select one PNR lot over another. Travel time in both the auto and transit networks is considered in the models, with data sourced from strategic transport networks and Google's transit feed specification. The study explores the trade-off relationship between the two networks and suggests that RRM models may better capture this compromise effect compared to traditional RUM models. Results show that RRM models, particularly the classical RRM, provide insights into transport choice behavior and offer improved predictive abilities compared to RUM models. Additionally, this study contributes to the understanding of RRM models by estimating them on revealed preference data.

The reference \cite{mao2020does} examines the use of RRM models in analyzing public preferences for air quality improvement policies, compared to traditional RUM models. The research finds that RRM models outperform RUM models in terms of goodness of fit and predictive accuracy. Approximately 55.9 \% of respondents' choices are influenced by regret. The study also reveals differences in preferences between regret-driven and utility-driven respondents: regret-driven individuals prioritize increased clean air days, reduced haze days, and lower mortality, while utility-driven individuals prioritize shortening policy delay. These findings can inform the design of socially optimal air quality policies and contribute to the development of environmental policies.

The reference \cite{iraganaboina2021evaluating} addresses the complex factors influencing route choice behavior in the context of advancing traffic management systems and the availability of diverse traffic information sources. While previous studies have focused on travel time and cost, this study examines the impact of information provision mechanisms on route preferences. Conducted in the Greater Orlando Region, USA, the research utilizes web-based stated preference surveys to gather data. The analysis employs RUM and RRM models to quantify the trade-offs among different attributes affecting route choice. Various variables, including trip characteristics, roadway features, and traffic information attributes, are incorporated into the models. The results offer insights into the trade-offs involved in travel information provision, aiding in the design of effective active traffic management systems. Additionally, the research customizes the computation for RRM models to account for variable interactions, enhancing the understanding of route choice behavior.

The reference \cite{wong2020revealed} investigates whether RRM models provide a better description of evacuation behavior compared to traditional RUM models. The researchers designed a revealed preference (RP) survey focusing on evacuation choices in the aftermath of the 2017 Southern California Wildfires. Contrary to their hypothesis, the data did not support the superiority of RRM over RUM for evacuation decisions. The lack of sufficient variation in attribute levels across alternatives made it challenging to distinguish between linear utility and non-linear regret. However, the study identified weak regret aversion for certain attributes and detected class-specific regret for route and mode choices using a mixed-decision rule latent class choice model. Despite the initial findings, the researchers suggest that RRM for evacuations could still be promising with further exploration. The study highlights methodological implications for RP studies with complex choice sets and limited attribute variability.

Additionally, DCM models can be considered as machine learning models. For instance, multinomial logit models correspond to well-known logistic regression models in ML.  Machine algorithms can be broadly categorized into one of the following three types: supervised, unsupervised, and agent-based learning, with DCMs falling into the supervised category. Paper \cite{boroujeni2024comprehensive} provides a comprehensive overview of the existing ML techniques along with their types, characteristics, and potential applications. Indeed, ML tools such as MNL and MMNL are used to predict customers' choices based on RUM theory \cite{hillel2021systematic, soleymani_forecasting_2024, boroujeni2024ic,  9814765, yousefpour2023unsupervised}. Reinforcement learning as an online learning method can also be leveraged to help RM \cite{jebellat2023reinforcement, ran2024development, amiri2023rank, amirimargavi2023lowrank, talebi2024leadership}.

Accurately predicting demand is a crucial aspect of revenue optimization, vital for effective supply planning and pricing strategies. Demand forecasting relies on various factors, such as customer diversity, pricing, competitor offerings, brand reputation, and product attributes. Traditionally, demand estimation employed a simple linear function, focusing solely on product price, which proved to be inadequate as it overlooked many influential factors. Additionally, the conventional assumption that customers faced only binary choices — either the product was available for purchase or not due to capacity constraints — did not reflect the complexities of real-world consumer behavior \cite{strauss2018review}. Recognizing these limitations, segmentation of customer offers emerged in the early 21st century, categorizing customers based on demographics like age, income, etc., and providing tailored choice sets and prices accordingly  \cite{talluri2004revenue,strauss2018review}. Confronted with limited capacity, stochastic dynamic programming models emerged to define product offer sets while accounting for capacity constraints over selling periods. Despite their efficacy, these models encountered difficulties in addressing large-scale problems \cite{strauss2018review}. Moreover, they overlooked customer individual-level heterogeneity, typically categorizing customers into broad groups.
The reference \cite{andersson1998passenger} made a pioneering contribution to Revenue Management by focusing on modeling passenger preferences using Discrete Choice Models (DCMs). They considered the heterogeneity of customers and estimated the probability of purchasing a product using DCMs, integrating it into their stochastic dynamic programming model. Following \cite{andersson1998passenger}, \cite{gallego2004managing} addressed a general question in RM related to consumer choice behavior (RUM-based). They explored the network revenue management problem with flexible products, considering scenarios where demand was either exogenously given or generated by a consumer choice model based on RUM-based DCM. Their findings indicated that the optimal value of the stochastic optimization problem closely approximated that of a deterministic linear program. The reference \cite{liu2008choice} delved into the implications of the concept introduced by \cite{gallego2004managing} in RM. They developed a Choice-based Deterministic Linear Programming (CDLP) model under consumer choice behavior, focusing on providing suitable offer sets, seats, and flights in the airline industry while adhering to capacity constraints over time. Their study demonstrated that the CDLP, inspired by \cite{gallego2004managing}, serves as a deterministic approximation of the stochastic formulation, with its solution being an asymptotic optimal solution for the original stochastic problem influenced by the stochasticity of customer choice probability. 
The reference \cite{bront2009column} expanded on the CDLP introduced by \cite{liu2008choice}. They addressed the issue of overlapping offer sets in the airline industry, where airlines typically present affordable morning tickets in one set and pricier afternoon tickets in another set. This extension considered combinations such as expensive and morning tickets, providing insights into more realistic scenarios. They applied the column generation method to tackle the problem's dimensionality issue. Particularly, meta-heauristc methods seem to be a promising approach in approximating high-quality solutions for large-sized problems (e.g., see \cite{Haeri_Boroujeni_2023, mehrabi2023efficient, boroujeni2021data, 9721496, mehrabi2021application}). The reference \cite{azadeh2015impact} introduced a modified CDLP model, solving it on a synthetic dataset to obtain optimal offer sets and an upper bound on revenue. Additionally, they suggested a non-parametric model for estimating choice probabilities, utilizing the synthetic data employed in the CDLP.

Thus far, we have elucidated the role of certain OR techniques, such as Stochastic programming and Linear programming, in RM. However, despite some models considering individual customer heterogeneity and RUM-based DCMs, they have typically assumed a fixed price for all customers or customer segments. This overlooks the influence of price dynamics on customer choice and demand, as demand is directly influenced by customer decisions. To address this limitation, the pricing of products should be treated as a variable rather than a fixed parameter in the models. This challenge was first tackled by \cite{bierlaire2016demand}, who introduced product prices as variables alongside considerations of individual customer heterogeneity and capacity constraints. Consequently, conventional RUM-based DCMs, which previously predicted the probability of a customer purchasing a product at a fixed price, transformed into nonlinear and non-convex demand functions due to the inclusion of price as a variable.
To mitigate this complexity, \cite{bierlaire2016demand} linearized the RUM-based DCMs using MILP models, thereby pioneering the application of RUM-based MILP models in RM. Subsequently, \cite{pacheco2017integrating} extended this approach by linearizing a mixed multinomial logit DCM, which considers unobserved customer utility correlated with customer attributes, through MILP techniques. Their objective was to maximize revenue by simultaneously optimizing supply (price and capacity) and demand (price, utility levels, customer satisfaction) factors. The reference \cite{paneque2019passenger} delved into an economics-oriented approach, emphasizing passenger satisfaction maximization by integrating DCMs into MILP models. The objective was to enhance the public transportation network by utilizing toll revenue from highways, ensuring a balance between economic considerations and passenger contentment. To date, no RRM-based MILP model has been applied for in RM.

In conclusion, DCMs are state-of-the-art for estimating dis-aggregate demand, showing promise in RM. The literature demonstrates integrating DCMs into optimization models for better supply-demand interactions. However, a gap exists between DCMs derived from theories like RRM and optimization formulations. 
Finally, \cite{pacheco2016new, pacheco2017integrating, PACHECOPANEQUE202126, PACHECOPANEQUE2022105985} developed a RUM-based DCM framework where customers are utility maximizers. We extended the framework for RRM-based DCM introduced by \cite{chorus2010new}, where the customers are regret minimizers.

\section{Modelling}
\label{sec3}
\subsection{Problem definition}

RRM posits that individuals select alternatives to minimize the regret they may encounter from foregoing other options. In other words, they conduct pairwise comparisons of the regret associated with each alternative or product, including the option of not choosing any alternative (opt-out), and ultimately choose the alternative expected to result in the least regret. Randomness in the context of RRM arises from researchers' inability to observe all the regret a customer experiences. To elucidate the RRM mechanism, consider the example of passenger travel mode choice provided by \cite{chorus2008random}. In this scenario, three alternatives are considered: car, train, and bus, each characterized by attributes $k^1$, $k^2$, and $k^3$, representing time, cost, and alternative-specific dummy variables of travel such as comfort level, reliability, accessibility, or any other characteristic that distinguishes one mode of transportation from another. The regret for each alternative is calculated as follows, and a rational customer would opt for the alternative with the minimum regret:

\begin{align*}
&R_i = \sum_{j\neq i} R_{ij}\ \text{and}\ R_{ij} = \sum_{k\in \{k^1,k^2,k^3\}}R_{ijk},\ \text{and}\ R_{ijk} = \max\{v_{ok}, \beta_k (k_j - k_i) + v_{k} \}&\\
\end{align*}

Where $\beta$ parameters, representing consumers' tastes for each attribute $k$, are estimated from historical data, and $k_i$ represents the $k^{th}$ attribute of alternative $i$. 
Regret, partially known to researchers, is modified by adding an error term, drawn from an independent and identically distributed (iid) extreme value type I distribution such as Gumbel distribution, to obtain the overall regret and $v_{ok}$ and $v_{k}$ stand for unobserved regret in attribute perceptions \cite{chorus2010new}. Thus, given two alternatives $i$ and $j$, we first compute the regret of choosing alternative $i$ attribute-wise, i.e., $R_{ijk}$, and then we obtain the sum of attribute-level regret, i.e., $R_{ij}$. Finally, by summing all the regret caused by selecting alternative $i$ given all other alternatives, we obtain the total anticipated regret caused by selecting alternative $i$, i.e., $R_i$. We will compute these values in MILP for each customer and all available alternatives for each customer.
Additionally, let $I = |C|$ represent the number of alternatives in the universal choice set $C$,  The set $N = \{1,...,n\}$ represents customers. For each customer $n$, the choice set $C_n$ is considered, with each product $i$ available in quantity $c_i$. Note that $i=0$ shows the opt-out alternative. 
The supplier aims to maximize revenue by making decisions on whether a product $i$ should be offered to customer $n$ and what price of alternative $i$ to be offered to customer $n$, denoted as $p_{in}$. 
Customers are served sequentially based on a priority list, determined by factors such as customer loyalty, arrival order, customer segments, anticipated order quantities, etc. Such an assumption is standard (e.g., see \cite{pacheco2016new, PACHECOPANEQUE202126}). As it was shown by \cite{pacheco2016new, paneque2017integrating}, when solving this problem directly using DCMs, highly non-linear and non-convex problems arise. Therefore, to address this challenge, we employ MILP models. This paper explores two cases, the first assumes the supplier has no capacity limitation on supplies, while the second considers a scenario where there is a capacity limitation on the supply.

\subsection{Model formulation}

We have summarized the notations used in the model in the following.

\subsubsection{Sets}
\begin{align*}
& C && \text{Set of all alternatives. $C_n$ is the set of all available alternatives for customer $n$} \\
& L_{in}  &&\text{Set of price levels for each alternative $i$ and customer $n$ indexed by $I$}\\ 
& N  &&\text{Set of customers, i.e., N= \{1,...,n\}} \\
& Ra  &&\text{Set of random draws, where each draw is indexed by  $r$ and $Ra = \{1,..,r\}$}
\end{align*}
Note that $C_n$ is of significance when the model is capacitated, as prior customers to customer 
$n$, for instance, can purchase all available quantity of an alternative. We assume that there is an order list from which the supplier satisfies the orders of the customers accordingly.To give an example of \( L_{in} \), suppose \( L_{in} = \{1,2,3\} \). Additionally, there is a fixed parameter for pricing, \( pm \), let us assume \( pm = 0.5 \). Then, the price levels of alternative \( i \) for customer \( n \) will be \( 0.5, 1, \) and \( 1.5 \), which are obtained by element-wise multiplication of \( pm \) by the elements of \( L_{in} \). \( pm \) allows us to adjust the price of alternatives for customers.
Finally, $Ra$, represents the set of scenarios or draws. This accounts for the stochastic behavior of customers, and we solve our models for several draws to obtain robust results.
\subsubsection{Input Parameters}

\begin{align*}
& c_i && \text{Available quantity of alternative } i \text{ for sale, maximum capacity of the alternative } i \\
& lp_{in}, mp_{in} && \text{Lower and upper bounds on the price, } lp_{0n} = mp_{0n} = 0, \text{i.e., the price of opt-out} \\ 
& && \text{alternative is zero.} \\
& pm && \text{The price multiplier is used to attain desired price levels. By adjusting the value }\\ 
& && \text{ of the multiplier, one can precisely set the price levels and ranges. Increasing or}\\
& && \text{decreasing the multiplier value allows for fine-tuning of prices.}\\
& ER_{ijnr}  && \text{Sum of attribute level regret associated with alternative $i$ when compared} \\ 
&     && \text{ with alternative $j$ for customer $n$ at scenario $r$. Note that:}  \\
& && ER_{ijnr} = \sum_x\sum_{k\in \{k_1,..., k_a\}} max\{v_{onxr}, \beta_{nx} (x_{kj}-x_{ki})+v_{nxkr}\}\\
& && \text{where $v_{onxr},v_{nxkr}$ are draws of Gumbel$(0,1)$ for customer $n$, and $k^{th}$ level of } \\
& && \text{attribute $x$ at draw $k$. Note that $x_{kj}$ represents the $k^{th}$ level of alternative $j$'s} \\ 
& &&\text{attribute $x$. Note that the price attribute has been excluded from $ER_{ijnr}$}\\
& && \text{since it is a variable in our models.} \\
& mm_{ijnr} && \text{Upper bounds on variables $RR_{ijnr}$ equivalent to $\max\{ v_{onr}, \beta_{np}(lp_{jn} - mp_{in}) + v_{nr} \}\}$} \\
 & ll_{ijnr} && \text{ Lower bounds on variables $RR_{ijnr}$ equivalent to $ \min\{ v_{onr}, \beta_{np}(mp_{jn} - lp_{in}) + v_{nr} \}\ $} \\
 & M_{ijnr} && \text{Big-M notion in MILP models, equivalent to } mm_{ijnr} - ll_{ijnr}  \\
 & m_{inr} && \text{Upper bounds on variables $R_{inr}$, which is } \text{$\sum_{j\neq i} \max\{ v_{onr}, \beta_{np}(lp_{jn} - mp_{in}) + v_{nr} \} + ER_{ijnr}$}\\
 & l_{inr} && \text{Lower bounds on variables $R_{inr}$, which is $\sum_{j\neq i} \min\{ v_{onr}, \beta_{np}(mp_{jn} - lp_{in} ) + v_{nr} \} + ER_{ijnr}$}\\
 & m_{nr}  &&   \text{Upper bounds on variables } R_{nr}, \textbf{equivalent to} \max_{i\in C_n} \{ m_{inr} \}\ \\
 &  l_{nr}   && \text{Lower bounds on variables $R_{nr}$ which is equivalent to $ \min_{i\in C_n} \{ l_{inr} \}$}\\
 & M_{nr} && \text{Big-M notion, equivalent to $m_{nr} - l_{nr}$} \\
 & M &&  \text{Big-M notion, equivalent to $ \max_{\forall i,j\in C_n, i\neq j, n,r }\{ |ll_{ijnr} + ER_{ijnr}|, |mm_{ijnr} + ER_{ijnr}|\}$} \\
& \beta_{np} && \text{Customer  $n^{th}$ taste for price attribute of alternatives.} \\
 & \beta_{nm} && \text{Customer $n^{th}$ taste for attribute $m$ of alternatives.}
\end{align*}
Note that all upper bound and lower bound parameters aid in potentially solving the MILP efficiently and quickly. However, their inclusion in the model is not mandatory. Specifically, setting upper bounds on Big-M, which assists in activating and deactivating the constraints of the MILP models, usually enhances solution time. Regarding \( \beta_{np} \), it indicates customer \( n \)'s preference for the price attribute. For instance, if \( \beta_{np} = -1 \), and the customer purchases the alternative, it implies that the customer is willing to pay one dollar for each unit of the alternative.

\subsubsection{Variables}
\begin{align*}
& y_{in} \in \{0,1\} && \text{Supplier’s decision when alternative $i$ is in the customer $n$’s choice set} \\
& y_{inr} \in \{0,1\} && \text{Availability variables, $y_{inr} = 1 $ if alternative $i$ is available to customer $n$ at scenario $r$} \\
& RR_{ijnr} && \text{Regret of customer $n$ for the price attribute of alternative $i$ when compared} \\ & && \text{ with that of alternative $j$} \\
& R_{inr} && \text{Overall regret of alternative } i \\
& R_{nr} && \text{Minimum regret among all available alternatives, i.e., $\min_{i\in C_n}R_{inr}$} \\
& b_{ijnr}\in \{0,1\} && \text{Binary variables to capture the $R_{inr}$}\\
& z_{ijnr}  && \text{Discounted regret, $z_{ijnr} \geq 0$ and $z_{ijnr} = RR_{ijnr} + ER_{ijnr} $ when $y_{jnr} = 1$ and 0 otherwise}\\
& w_{inr} \in \{0,1\} && \text{Choice of the customer, } w_{inr} = 1 \text{ if service } i \text{ is chosen by customer } n \text{ at scenario } r. \\
& \lambda_{inl} \in \{0,1\} && \text{Binary variables to identify price levels}. \\
& \alpha_{inrl} \in \{0,1\} && \text{Linearization variables }.
\end{align*}
Among the variables, \( y_{in} \) relates to the direct decision of the supplier to grant access to alternative \( i \) for customer \( n \) or not. However, \( y_{inr} \) relates to the availability of alternative \( i \) for customer \( n \) in scenario \( r \). In other words, it not only considers the supplier's decision regarding the availability of alternative \( i \) for customer \( n \), but also takes into account the available capacity of alternative \( i \) through the constraints of the MILP model. 
Variables \( \alpha_{inrl} \) are technical variables used to linearize the objective function of the MILP models presented in this paper. When \( \alpha_{inrl} = 1 \), it indicates that customer \( n \) purchases alternative \( i \) in scenario \( r \) at price level \( l \); otherwise, it is equal to zero. Similarly, variables \( b_{ijnr} \) are related to techniques in MILP models to activate or deactivate constraints. In our paper, these variables assist in computing the anticipated regret of purchasing alternative \( i \) by customer \( n \) in scenario \( r \). Variables \( z_{ijnr} \) assist us in computing the discounted regret of choosing alternative \( i \) over alternative \( j \) by customer \( n \) in scenario \( r \). We refer to it as "discounted" following Pacheco et al. (2016) \cite{pacheco2016new}, as the purchasing decision has not yet been finalized by the customer, and we are anticipating the regret. The remaining variables pertain to the calculation of regret and the selection of the final choice. 

\subsubsection{Uncapacitated supply MILP model}

The relationships between variables are captured through the MILP models. In this part, we first present our models, and then we will explain the constraints and relationships between variables.

\begin{align}
&\max \frac{1}{|Ra|}\sum_n\sum_{i>0,i\in C_n}\sum_r(lp_{in}w_{inr}+\sum_l l\alpha_{inrl}pm)\label{eqn:obj}\\
&s.t.\nonumber \\
&v_{onr}\leq RR_{ijnr},&\forall i,j\in C_n, i\neq j, n, r,\label{r1}\\
&RR_{ijnr}\leq v_{onr}+M_{ijnr}b_{ijnr},&\forall i,j\in C_n, i\neq j, n, r,\label{r2}\\
&\beta_{np}(lp_{jn}+\sum_l pm \lambda_{jnl}l-lp_{in}-\sum_l pm\lambda_{inl}l)\nonumber \\
&+v_{nr}\leq RR_{ijnr},&\forall i,j\in C_n, i\neq j, n, r,\label{r3}\\
&RR_{ijnr}\leq \beta_{np}(lp_{jn}+\sum_l l\lambda_{jnl}pm-lp_{in}-\sum_l l\lambda_{inl}pm)\nonumber \\
&+v_{nr}+M_{ijnr}(1-b_{ijnr}),&\forall i,j\in C_n, i\neq j, n, r,\label{r4}\\
&R_{inr}=\sum_{j\neq i, j\in C_n}RR_{ijnr}+ER_{ijnr},&\forall i\in C_n, n, r,\label{tot}\\
&R_{nr}\leq R_{inr}, &\forall i\in C_n, n, r\label{c1}\\
&R_{inr}-M_{nr}(1-w_{inr})\leq R_{nr},&\forall i\in C_n, n, r\label{c2}\\
&\alpha_{inrl}\leq \lambda_{inl},&\forall i\in C_n,i\neq 0,n,r,l,\label{obj1}\\
&\alpha_{inrl}\leq w_{inr},&\forall i\in C_n,i\neq 0,n,r,l,\label{obj2}\\
&\sum_l \lambda_{inl}\leq 1,&\forall i\in C_n,i \neq 0, n,\label{obj3}\\
&lp_{in}+\sum_l l\lambda_{inl}pm\leq mp_{in}&\forall i\in C_n, i\neq  0, n\label{last}
\end{align}

There are \(|Ra|\) scenarios for each customer in the objective function \ref{eqn:obj}. This is due to the stochastic nature of consumer behavior. Averaging over several scenarios will result in more robust solutions. The \(lp_{in}\) variable represents a lower bound on the price levels, i.e., the minimum price at which alternative \(i\) will be offered to customer \(n\). Additionally, exactly one price level should be selected in each scenario \(r\), which maximizes revenue while respecting the constraints of the model, i.e., minimizing customer regret. 
If alternative \(i\) is chosen (i.e., \(w_{inr} = 1\)) for the customer \(n\) in scenario \(r\), the corresponding price level (\(\lambda_{inl} = 1\)) should be selected in such a way that \(lp_{in} + \sum_l l\alpha_{inrl}pm\) contributes most positively to the objective function, while still respecting the constraints of the model, as it is a maximization problem.
We define \(pm\) as a constant price multiplier, and integers \(l\) are used to characterize price levels. The opt-out alternative is excluded from the objective function, (i.e., $i>0$) as it does not affect revenue (i.e., the price of the opt-out alternative is zero).

To calculate the anticipated regret arising from the comparison of the price attribute of two alternatives, we introduce constraints \ref{r1},\ref{r2}, \ref{r3}, and \ref{r4}. In the uncapacitated version, the above MILP model, all alternatives are available for each customer. An error term \( v_{nr} \) has been added to the price-raised regret between alternatives \( i \) and \( j \), \( RR_{ijnr} \). Finally, we compare this anticipated regret with that of unobserved heterogeneity, \( v_{onr} \), to determine the maximum anticipated regret. In more details, 
\begin{align*}
RR_{ijnr} = \max\{v_{onr}, \beta_{np}(lp_{jn} + \sum_{l} l\lambda_{jnl}pm - lp_{in} - \sum_{l} l\lambda_{inl}pm) + v_{nr}\}, \quad \forall i, j \in C_n, i \neq j, n, r,
\end{align*}
where $lp_{jn} + \sum_{l} l\lambda_{jnl}pm$ characterizes the price of the alternative $j$ for customer $n$. 
Since \( RR_{ijnr} \) is the maximum of two terms, it can be explicitly asserted that:
\begin{align*}
& v_{onr} \leq RR_{ijnr}, & \forall i, j \in C_n, i \neq j, n, r, \\
& \beta_{np}(lp_{jn} + \sum_{l} lp_{jn} \lambda_{jnl} - lp_{in} - \sum_{l} lp_{in} \lambda_{inl}) + v_{nr} \leq RR_{ijnr}, \quad & i, j \in C_n, i \neq j, n, r 
\end{align*}
However, \( RR_{ijnr} \) is exactly equal to the maximum of these two terms, which is ensured by constraints \ref{r2} and \ref{r4} as follows:\\
If \( b_{ijnr} = 1, \forall i, j \in C_n, i \neq j, n, r \), it is obvious that:
\begin{align*}
&RR_{ijnr} \leq v_{onr} + M_{ijnr}, &\forall i, j \in C_n, i \neq j, n, r 
\end{align*}
These constraints will be inactive since:
\begin{align*}
&ll_{ijnr} \leq v_{onr}, &\forall i, j \in C_n, i \neq j, n, r \\
& mm_{ijnr} - ll_{ijnr} = M_{ijnr}&\\
&RR_{ijnr} \leq mm_{ijnr} \leq v_{onr} + mm_{ijnr} - ll_{ijnr} = v_{onr} + M_{ijnr}, & \forall i, j \in C_n, i \neq j, n, r 
\end{align*}
meaning that constraints \ref{r2} always hold if \( b_{ijnr} = 1, \forall i, j \in C_n, i \neq j, n, r \).
At the same time, due to constraints \ref{r4}, when \( b_{ijnr} = 1, \forall i, j \in C_n, i \neq j, n, r \),
\begin{align*}
RR_{ijnr} \leq \beta_{np}(lp_{jn} + \sum_{l} lp_{jn} \lambda_{jnl} - lp_{in} - \sum_{l} lp_{in} \lambda_{inl}) + v_{nr}, \quad \forall i, j \in C_n, i \neq j, n, r
\end{align*}
Due to constraints \ref{r3}, we also have: 
\begin{align*}
 & \beta_{np}(lp_{jn} + \sum_{l} lp_{jn} \lambda_{jnl} - lp_{in} - \sum_{l} lp_{in} \lambda_{inl}) + v_{nr} \leq RR_{ijnr}, & \forall i, j \in C_n, i \neq j, n, r
\end{align*}
This implies that 
\begin{align*}
 & \beta_{np}(lp_{jn} + \sum_{l} lp_{jn} \lambda_{jnl} - lp_{in} - \sum_{l} lp_{in} \lambda_{inl}) + v_{nr} = RR_{ijnr}, & \forall i, j \in C_n, i \neq j, n, r
\end{align*}
when 
\begin{align*}
 & v_{onr} \leq \beta_{np}(lp_{jn} + \sum_{l} lp_{jn} \lambda_{jnl} - lp_{in} - \sum_{l} lp_{in} \lambda_{inl}) + v_{nr} 
\end{align*}
With the same rationale, when $b_{ijnr}=0$, we have $v_{onr}= RR_{ijnr}$. Indeed, these sets of constraints help us to linearize 
\begin{align*}
RR_{ijnr} = \max \{v_{onr}, \beta_{np}(lp_{jn} + \sum_{l} lp_{jn} \lambda_{jnl} - lp_{in} - \sum_{l} lp_{in} \lambda_{inl}) + v_{nr}\}
\end{align*}

Considering the RUM-based MILP, we compare alternatives pairwise to calculate regret. We compute the regret obtained from different price levels for each alternative pairwise, i.e., in comparison to other alternatives. Then, the alternative rendering the overall minimum regret will be chosen.
The parameters \( ER_{ijnr} \) represent the sum of the anticipated regret of attributes (other than the price attribute) for alternatives, and constraints \ref{tot} are used to obtain the overall anticipated regret of an alternative by summing \( ER_{ijnr} \) and the regret related to the price attribute, i.e., \( RR_{ijnr} \). Note that \( ER_{ijnr} \) are parameters of the model and they have been calculated beforehand.

It is very unlikely to have two alternatives with the minimum regret. Without loss of generality, we assume that ties are not allowed. Constraints \ref{tot}, \ref{c1} and \ref{c2} in our model, capture the alternative with the lowest regret as follows:
\begin{align}
&R_{inr} = \sum_{j\neq i, j\in C_n} RR_{ijnr} + ER_{ijnr},& \forall i \in C_n, n, r, \label{eq:Rinr} \\
&R_{nr} \leq R_{inr}, & \forall i \in C_n, n, r, \label{eq:Rnr} \\
&R_{inr} - M_{nr}(1 - w_{inr}) \leq R_{nr}, &\forall i \in C_n, n, r \label{eq:Rinr_M}
\end{align}
Variables \( R_{nr} \) capture the minimum regret among all alternatives in the choice set of customer \( n \) at scenario \( r \). The overall regret of alternative \( i \) for the customer \( n \) at scenario \( r \) is measured by variable \( R_{inr} \) through constraint \eqref{eq:Rinr}. 
 Let us assume that \( i \) is the alternative with the minimum regret at scenario \( r \) for customer \( n \) in her choice set. So, \( w_{inr} = 1 \), i.e., the customer $n$'s choice at scenario $r$ is alternative $i$ as it brings about the least anticipated regret among alternatives. Due to constraints \eqref{eq:Rnr} and \eqref{eq:Rinr_M}, by substituting $w_{inr} =1$, we have
\[
R_{nr} \leq R_{inr} \quad \text{and} \quad R_{inr} \leq R_{nr}
\]
Which can be inferred that,
\[
R_{nr} = R_{inr}
\]
There is not any other alternative such as \( j \) that \( w_{jnr} = 1 \). Otherwise, \( R_{inr} \leq R_{nr} \) and \( R_{jnr} \leq R_{nr} \). Moreover, we have \( R_{nr} \leq R_{inr} \) and \( R_{nr} \leq R_{jnr} \). We can conclude that \( R_{nr} = R_{inr} = R_{jnr} \). This is in contradiction with the assumption that ties are not allowed. Hence, there exists exactly one alternative such as \( i \) for which \( w_{inr} = 1 \) and \( R_{nr} = R_{inr} \). Hence, \( w_{jnr} = 0, \forall j \neq i, j \in C_n \). It can be checked that constraints \eqref{eq:Rinr}, \eqref{eq:Rnr}, and \eqref{eq:Rinr_M} are inactive when \( w_{jnr} = 0 \).
So far, we have shown how the alternative with the minimum anticipated regret will be selected by a customer in different scenarios. However, note that we have assumed that the prices of the alternatives have been chosen by the other constraints. In other words, there are intricate relationships between pricing and the choice of the customers. We have just covered how the anticipated regret of alternatives has been calculated. Part of this regret comes from the price of the product, which affects constraints \eqref{r1}, \eqref{r2}, \eqref{r3}, and \eqref{r4} by influencing the \(RR_{ijnr}\) through variables $\lambda_{inl}$ in constraints \ref{r2}. For each alternative $i$ and customer $n$, this binary variable will take a value of 1 for only one price level, i.e., $\sum_l \lambda_{inl}=1$, as indicated by constraints \eqref{obj3}. Due to these constraints and since it is a maximization problem, it is certain that $\lambda_{inl}$ and $w_{inr}$ will each be equal to exactly one for an alternative $i$ at scenario $r$ and a price level $l$, respecting other constraints of the model. Thus, constraints \eqref{obj1} and \eqref{obj2} will take the form $\alpha_{inrl}\leq 1$. Again, since it is a maximization problem, $\alpha_{inrl} = 1$ so that it has the most possible positive contribution to the objective function. This is how constraints help us compute the regrets and select the alternative with the lowest anticipated regret while maximizing the revenue and pricing the products.
Note that the constraints \eqref{obj3} help us establish an upper bound for the price of a specific alternative $i$ and customer $n$ by not allowing the chosen price to exceed the threshold, i.e., $mp_{in}$. This is useful when categorizing our customers; for example, the product price may be cheaper for students.

The remaining constraints are bound constraints, where the bounds were defined in the notation part. 

\begin{align}
&RR_{ijnr} \leq mm_{ijnr} &\forall i, j \in C_n, i \neq j, n,r , \\
&ll_{ijnr} < RR_{ijnr} &\forall i, j \in C_n,i \neq j,  n,r,  \\
&R_{inr} \leq m_{inr} &\forall i, j \in C_n, n,r , \\
&l_{inr} < R_{inr} &\forall i \in C_n, n,r , \\
&R_{inr}, RR_{ijnr}\in \mathbb{R} &\forall i, j \in C_n, i \neq j, n,r. 
\end{align}

These bounds do not affect the quality of the solutions; however, these constraints typically enhance the efficiency of the model-solving process.

\subsubsection{Capacitated supply MILP model}
A capacitated version of the previously presented model has been depicted in this part. The term "capacitated" is used to emphasize that the available quantity of an alternative is limited.

\begin{align}
&\max \frac{1}{R}\sum_n\sum_{i\neq 0,i\in C_n}\sum_r(lp_{in}w_{inr}+\sum_l l\alpha_{inrl}pm)\label{eqn:91}& \\
&s.t.&\nonumber \\
&y_{in}=0,&\forall i \notin C_n,n,r,\label{eqn:92}\\
&y_{inr}\leq y_{in},&\forall i,n,r,\label{eqn:93}\\
&v_{onr}\leq RR_{ijnr},&\forall i,j\in C_n, i\neq j, n, r,\label{eqn:94}\\
&RR_{ijnr}\leq v_{onr}+M_{ijnr}b_{ijnr},&\forall i,j\in C_n, i\neq j, n, r,\label{eqn:95}\\
&\beta_{np}(lp_{jn}+\sum_l  l\lambda_{jnl}pm-lp_{in}-\sum_l l\lambda_{inl}pm)\nonumber\\
&+v_{nr}\leq RR_{ijnr},&\forall i,j\in C_n, i\neq j, n, r,\label{eqn:96}\\
&RR_{ijnr}\leq \beta_{np}(lp_{jn}+\sum_l l\lambda_{jnl}pm-lp_{in}\nonumber \\
&-\sum_l  l\lambda_{inl}pm)+v_{nr}+M_{ijnr}(1-b_{ijnr}),&\forall i,j\in C_n, i\neq j, n, r,\label{eqn:97}\\
&z_{ijnr}\leq RR_{ijnr}+ER_{ijnr}&\forall i,j\in C_n, i\neq j, n, r,\label{eqn:98}\\
&RR_{ijnr}+ER_{ijnr}-M(1-y_{jnr})\leq z_{ijnr}&\forall i,j\in C_n, i\neq j, n, r,\label{eqn:99}\\
&z_{ijnr}\leq My_{jnr}&\forall i,j\in C_n, i\neq j, n, r,\label{eqn:100}\\
&-My_{jnr}\leq z_{ijnr}&\forall i,j\in C_n, i\neq j, n, r,\label{eqn:101}\\
&R_{inr}= \sum_{j\neq i,j\in C_n} z_{ijnr}, &\forall i\in C_n, n, r,\label{eqn:102}\\
&w_{inr}\leq y_{inr},&\forall i\in C_n, n, r,\label{eqn:34}\\
&R_{nr}\leq R_{inr}+M_{nr}(1-y_{inr}), &\forall i\in C_n, n, r\label{eqn:105}\\
&R_{inr}-M_{nr}(1-w_{inr})\leq R_{nr},&\forall i\in C_n, n, r\label{eqn:106}\\
&y_{i(n+1)r}\leq y_{inr},&\forall  i\in C_n,i\neq 0,n<N,r,\label{eqn:107}\\
&\sum_{m=1}^{n-1}w_{imr}\leq (c_i-1)y_{inr}+(|N|-1)(1-y_{inr}),&\forall i\in C_n,i\neq 0, n>c_i,r,\label{eqn:108}\\
&c_i(y_{in}-y_{inr})\leq \sum_{m=1}^{n-1}w_{imr},&\forall  i\in C_n,i\neq 0,n,r,\label{eqn:109}\\
&\alpha_{inrl}\leq \lambda_{inl},&\forall  i\in C_n,i\neq 0,n,r,l,\label{eqn:110}\\
&\alpha_{inrl}\leq w_{inr},&\forall  i\in C_n,i\neq 0,n,r,l,\label{eqn:111}\\
&\sum_l \lambda_{inl}\leq 1,&\forall  i\in C_n,i\neq 0,n,\label{eqn:112}\\
&lp_{in}+\sum_l pm\lambda_{inl}l\leq mp_{in}&\forall  i\in C_n,i\neq 0, n\label{eqn:113}
\end{align}

Constraints \ref{eqn:92}, \ref{eqn:93}, \ref{eqn:34} are incorporated into the model to ensure that once an alternative is not in the choice set of a customer, it will never be allocated to that customer. When $y_{in} = 0$, variables $y_{inr}$ are imposed to be zero \eqref{eqn:93} and consequently, $w_{inr} = 0$ due to constraints \eqref{eqn:34}.

Constraints \ref{eqn:107}, \ref{eqn:108}, and \ref{eqn:109} ensure that the capacity of the alternatives is respected. Here is how it works: Since in this version of the model, we assume that all the alternatives are available to all customers unless capacity is limited if $y_{inr} = 0$, it means that there is not enough quantity of alternative $i$ for the customer $n$ at scenario $r$. Thus, it will not be available for the customer $n+1$, and consequently, $y_{i(n+1)r} = 0$ due to constraints \ref{eqn:107}. Then, constraints \ref{eqn:108} and \ref{eqn:109} will imply that there is not enough quantity for the alternative $i$, since constraint \ref{eqn:108} will be $\sum_{m=1}^{n-1}w_{imr} \leq n-1$, meaning that at most $n-1$ unit of alternative $i$ has been allocated to earlier customers (customers 1 to $n-1$), and constraint \ref{eqn:109} will be $c_i \leq \sum_{m=1}^{n-1}w_{imr}$, meaning that all $c_i$ quantity has been allocated to customers before customer $n$. If $y_{inr} = 1$, constraints \ref{eqn:107} imply that the customer $n+1$ may still purchase alternative $i$. Also, constraints \ref{eqn:108} will imply that $\sum_{m=1}^{n-1}w_{imr} \leq c_i -1$, i.e., at most $c_i-1$ unit of alternative $i$ has been allocated to customers earlier than customer $n$, and still, there is some capacity of alternative $i$ to be allocated to customer $n$.

Constraints \ref{eqn:94}, \ref{eqn:95}, \ref{eqn:96} and \ref{eqn:97} have been explained previously. However, the key difference between the capacitated model and the uncapacitated model is that there might not be an available quantity of an alternative in the capacitated model. Hence, the anticipated regret should be computed based on the available alternatives for the customer.
Constraints associated with discounted regret, \ref{eqn:98},\ref{eqn:99}, \ref{eqn:100}, \ref{eqn:101}, and \ref{eqn:102} help find the available alternatives with the least anticipated regret in the choice set of a customer. Therefore, we propose $z_{ijnr}$, discounted regret variables, $z_{ijnr} = R_{ijnr}$ while $y_{jnr} = 1$ and $z_{ijnr} = 0$ when $y_{jnr} = 0$. In other words, when the competing alternative is not available due to capacity limitations, there is no regret, and the anticipated regret should be computed among available alternatives. The idea is similar to that of RUM-based MILP developed by \cite{pacheco2016new}. However, in our model, since the regret is a result of the pairwise comparison of alternatives, it is important to check if both alternatives are available. Constraints  \ref{eqn:98},\ref{eqn:99}, and \ref{eqn:100} work as follows:\\
 If $y_{jnr} = 1$, then
\begin{align*}
&z_{ijnr} \leq RR_{ijnr} + ER_{ijnr}  \\
&RR_{ijnr} + ER_{ijnr} \leq z_{ijnr} 
\end{align*}
implying that
\begin{align*}
z_{ijnr} = RR_{ijnr} + ER_{ijnr}.
\end{align*}
Otherwise,
\begin{align*}
&z_{ijnr} \leq RR_{ijnr} + ER_{ijnr}\\
&RR_{ijnr} + ER_{ijnr} - M \leq z_{ijnr} 
\end{align*}
We know that $M$ is big enough, so
\begin{align*}
RR_{ijnr} + ER_{ijnr} - M \leq 0.
\end{align*}
Thus, $z_{ijnr}$ can take any value in $[RR_{ijnr} + ER_{ijnr} - M, RR_{ijnr} + ER_{ijnr}]$. But, we want this variable to be zero. Constraints \ref{eqn:100} ensure that while $y_{jnr} = 0$, then $z_{ijnr} = 0$. They will be inactive once $y_{jnr} = 1$ since $0 \leq z_{ijnr} \leq M$ always holds. Consequently, to avoid the effect of the not-available alternative $j$ on the overall anticipated regret of the alternative $i$, we replace the constraint set
\begin{align*}
RR_{ijnr} = \sum_{j \neq i, j \in C_n} RR_{ijnr} + ER_{ijnr}, \quad \forall i \in C_n, n, r, 
\end{align*}
with the constraint set
\begin{align*}
RR_{ijnr} = \sum_{j \neq i, j \in C_n} z_{ijnr}, \quad \forall i \in C_n, n, r. 
\end{align*}
Constraints \ref{eqn:34}, \ref{eqn:105}, \ref{eqn:106} address the issue of picking up the alternative with the least regret while it has enough quantity to be allocated.

The rest of the constraints are the same as in the previous model. Finally, constraints \ref{eqn:92}, \ref{eqn:93}, \ref{eqn:107}, \ref{eqn:108}, \ref{eqn:109}, \ref{eqn:110}, \ref{eqn:111}, \ref{eqn:113} were taken from \cite{pacheco2016new}.

\section{Computational experiments}\label{sec4}

In our experiments, we included three alternatives in the choice set for each customer. One of them represents the no-buy option, while the other two alternatives are left unlabeled. Regret arises due to the price attribute and unobserved factors. We did not consider the socio-economic characteristics of the customers in this experiment, i.e., $ER_{ijnr} =0,\ \forall i,j,n,r$. Specifically, price is the only attribute for both alternatives. Customers are homogeneous in their preferences,
so $\beta_{np} = -1,\ \forall n.$ Also, $v_{onr}$ and $v_{nr}$ are positive random draws from $Gumbel(0,1)$ distribution. Price has six levels ranging from $1$ to $4.5$, with $pm = 0.5$. The bounds are set to be arbitrarily large enough numbers. As for the RUM-based MILP modeled by \cite{pacheco2017integrating}, the overall utility is solely associated with the price attribute and the unobserved utility, which is a random draw from $Gumbel(0,1)$.

Firstly, we solved both capacitated and uncapacitated RRM-based MILP models to optimality. We considered $r= 4$ draws for each customer in both models, i.e., we considered 4 scenarios for each customer.  

\subsection{Numerical results}

We employ IBM ILOG CPLEX Optimization solver to execute the models and assess their performance in solving the numerical experiments.

The Table \ref{t11} summarizes the performance of both capacitated and uncapacitated models. The maximum revenue achievable in our experiments could be obtained by multiplying the number of customers by the highest price level, i.e., $\$\ 4.5$. Thus, all the results in the table imply that once there is a sufficient quantity of alternatives for each customer, maximum revenue is attained. In other words, regret stemming from the price of the alternatives is negligible compared to regret when choosing the opt-out option. This could be partly because we have only included the regret arising from the price of the alternatives.
To verify if capacity constraints are active, consider the instance with 11 customers. In this case, there are five units of each alternative available. Consequently, only the opt-out choice is available for one of the customers, which is reflected in the gap between the uncapacitated and capacitated solutions.\\
It is also important to note that despite having the same number of customers and equal capacity for the alternatives (except in the case of 11 customers), the uncapacitated model was solved more quickly. This could be attributed to the fact that the capacitated version involves additional constraints to adhere to capacity restrictions, leading to delays in solving the model.

\begin{table}[H]
    \centering
    \caption{Performance of Capacitated and Uncapacitated RRM-based MILP Models \label{t11}}
    \begin{tabular}{ccccccc}
        \hline
        No. Customer  & \multicolumn{2}{c}{Capacity} & \multicolumn{2}{c}{Time (s)} & \multicolumn{2}{c}{Solution} \\
        & Alt. 2 & Alt. 3 & Cap & Uncap & Cap & Uncap \\
        \hline
        10 & 5 & 5 & 171 & 2.2 & 45 & 45 \\
        11 & 5 & 5 & 216 & 3.4 & 45 & 49.5 \\
        12 & 6 & 6 & 206 & 3.5 & 54 & 54 \\
        13 & 7 & 6 & 357 & 3.6 & 58.5 & 58.5 \\
        14 & 7 & 7 & 787 & 9.5 & 63 & 63 \\
        15 & 8 & 7 & 1320 & 5.28 & 67.5 & 67.5 \\
        \hline
    \end{tabular}
     
\end{table}

In Table \ref{fig:2}, the number of constraints and variables for both models shows that the capacitated model contains more variables and constraints. This is due to the introduction of capacity-related constraints with additional variables $z_{ijnr}$, $y_{in}$, and $y_{inr}$. The additional variables and constraints, stemming from the combinatorial nature of these models, generally require more time for the model to be solved, as described in Table \ref{fig:2}. Also, since we did not label the alternatives, and customers are identical, there could be a lot of symmetries in the model. This means that the solution method of the model distinguishes between the same solutions, thereby increasing solving time.
The branch and cut method was used by the CPLEX solver to find the solutions. Indeed, the 'nodes' column in Table \ref{fig:2} demonstrates the number of nodes on the search tree of the branch and cut algorithm. A large number of nodes for the capacitated version shows that the algorithm assessed the same solutions many times due to the symmetry of the model, most possibly because of the capacity-related constraints, as these numbers are fewer in the uncapacitated model. Additionally, 'iterations' show how many times the branch and cut algorithm was used to solve the model. Certainly, if the number of nodes is high, the number of iterations will be high as well, since each node inspected by the algorithm implies additional runs of the algorithm.\\
All in all, the addition of capacity-related constraints has increased the solution time, which was expected due to the combinatorial nature of the models.

\begin{table}[h]
    \centering
        \caption{Structure and Performance of Capacitated and Uncapacitated RRM-based MILP Models  \label{fig:2}}
    \begin{tabular}{cccccccccccc}
        \hline
        No. Customer  & \multicolumn{2}{c}{Constraints} & \multicolumn{2}{c}{Variables} & \multicolumn{2}{c}{Iterations} & \multicolumn{2}{c}{Nodes} \\
        & Cap & Uncap & Cap & Uncap & Cap & Uncap & Cap & Uncap \\
        \hline
        10 & 5412 & 2520 & 1910 & 1530 & 10291738 & 5471 & 175717 & 0 \\
        11 & 5991 & 2772 & 2112 & 1683 & 5466747 & 7312 & 67791 & 6355 \\
        12 & 6532 & 3024 & 2304 & 1836 & 8627555 & 8752 & 147674 & 0 \\
        13 & 7077 & 3276 & 2496 & 1989 & 6352511 & 8538 & 67486 & 0 \\
        14 & 7622 & 3528 & 2688 & 2142 & 12983867 & 343831 & 222918 & 0 \\
        15 & 8167 & 3780 & 2880 & 2295 & 18771065 & 13401 & 169710 & 0 \\
        \hline
    \end{tabular}

\end{table}

In Table \ref{fig:3}, we have shown the results of the experiment for the uncapacitated version of the RRM-based and RUM-based MILP models. These models represent different theories of customer decision-making: the RRM-based model and the RUM-based model taken from \cite{pacheco2016new}.\\
In the RRM-based model, which is rooted in behavioral economics, customers are assumed to make decisions based on minimizing regret. This means they choose options that minimize the disappointment or dissatisfaction they might feel after making a decision. On the other hand, the RUM-based model is based on traditional economic theory and assumes that customers make decisions to maximize their utility or satisfaction.\\
The results presented in Table \ref{fig:3} indicate that in the RRM-based model, higher price levels tend to be selected by customers. This suggests that customers in this model are willing to pay more to avoid regret. Conversely, in the RUM-based model, customers favor lower price levels. This indicates that customers in this model are more focused on maximizing their utility by seeking lower prices.\\
This discrepancy highlights the impact of different choice behavior theories on results, emphasizing the importance of utilizing these theories appropriately in demand and supply studies. Applying RUM-based models where customers are regret minimizers can lead to significant revenue losses in the real world, as evidenced by the observed differences in revenue (the 'Gap' column) between the two models.\\
Note that in our experiments, we have not assumed the true behavior of the customers, i.e., we have not pre-assumed that customers are regret minimizers or utility maximizers in order to obtain unbiased results. Our goal is to demonstrate the potential implications of different customer decision-making theories, such as RRM and RUM, on the outcomes of our models.\\
The observed 'Gap' in the solutions of RRM-based and RUM-based models underscores the importance of understanding and accurately modeling customer behavior. If customers are indeed regret minimizers, as suggested by the results of the RRM-based model, then implementing such models could benefit both suppliers and customers. Suppliers can optimize revenue by offering price levels that minimize customer regret, while customers benefit from choices that align with their decision-making preferences.\\
Conversely, applying an RRM-based model in scenarios where customers are utility maximizers, as indicated by the results of the RUM-based model, may not be favorable for customers. In such cases, customers may perceive the offerings as less attractive or may experience dissatisfaction due to higher prices, potentially leading to revenue losses for suppliers.\\ 
The fact that the number of customers affects the number of variables, thereby increasing the size of the search tree in the Branch and Cut algorithm, results in longer solving time. This effect, however, is more pronounced for the RRM-based MILP as it contains more variables than its RUM-based counterpart. However, neither model provides insight into the regret or utility levels of customers for different price levels. Such information can significantly help suppliers devise more effective marketing and production strategies. For instance, developing bi-objective MILPs enables the generation of Pareto frontiers, which depict the trade-offs between prices and revenue for various levels of customer regret or utility. This allows suppliers to make informed decisions that optimize both revenue and customer satisfaction.

\begin{table}[H]
    \centering
         \caption{Performance of Uncapacitated RRM and RUM-based MILP Models \label{fig:3}}
         \resizebox{\textwidth}{!}{
\setlength{\tabcolsep}{7pt}
\blt
    \begin{tabular}{ccccccc}
        \hline 
        No. Customers & Time, RUM (s) & Time, RRM (s) & Solution, RUM & Solution, RRM & Gap(\%) \\
        \hline
        10 & 0.08 & 2.2 & 7.125 & 45 & 84 \\
        11 & 0.02 & 3.4 & 21.125 & 49.5 & 57 \\
        12 & 0.08 & 3.5 & 10 & 54 & 81 \\
        13 & 0.06 & 3.6 & 8.25 & 58.5 & 85 \\
        14 & 0.06 & 9.5 & 7.25 & 63 & 88 \\
        15 & 0.08 & 5.28 & 15 & 67.5 & 77 \\
        \hline
    \end{tabular}
}
\end{table}

In Table \ref{fig:4}, the structural properties of RUM-based and RRM-based MILPS for several instances are presented. As the number of customers increases, both models exhibit an increase in the number of variables and constraints. However, the RRM-based MILP demonstrates a larger number of variables and constraints compared to its RUM-based counterpart. This disparity arises from the additional constraints and variables required by the RRM theory, which need to be linearized and integrated into the MILP model.\\
Furthermore, Table \ref{fig:4} reveals that the solver required significantly more iterations to solve the RRM-based MILP compared to its RUM-based counterpart. This indicates that for large-sized problems with numerous customers, the RRM-based model will require a substantial amount of time to achieve optimality. Consequently, operations research techniques such as metaheuristic methods and decomposition techniques should be developed so that the RRM-based MILP can obtain approximate or exact solutions in a timely manner. 
 
\begin{table}[!h]
    \centering
    \caption{Structure and Performance of Uncapacitated RRM and RUM-based MILP Models \label{fig:4}}
    \resizebox{\textwidth}{!}{
    \setlength{\tabcolsep}{7pt}
    \blt
    \begin{tabular}{cccccccc}
        \hline
        No. Customer & \multicolumn{2}{c}{Constraints} & \multicolumn{2}{c}{Variables} & \multicolumn{2}{c}{No. Iter} \\
        & RUM-based & RRM-based & RUM-based & RRM-based & RUM-based & RRM-based \\
        \hline
        10 & 1570 & 2520 & 1050 & 1530 & 26 & 5471 \\
        11 & 1727 & 2772 & 1155 & 1683 & 17 & 7312 \\
        12 & 1884 & 3024 & 1260 & 1836 & 23 & 8752 \\
        13 & 2041 & 3276 & 1365 & 1989 & 27 & 8538 \\
        14 & 2198 & 3528 & 1470 & 2142 & 12 & 343831 \\
        15 & 2355 & 3780 & 1575 & 2295 & 47 & 13401 \\
        \hline
    \end{tabular}
}
\end{table}

\section{Conclusion}\label{sec5}

In conclusion, cutting-edge Revenue Management (RM) techniques, adept at handling the intricate dynamics of supply and demand, leverage Discrete Choice Models (DCMs) as superior demand estimators rooted in consumer purchasing behavior. These models are seamlessly integrated into Mixed-Integer Linear Programming (MILP) frameworks, drawing inspiration from Random Utility Maximization (RUM) theory and the more recent Random Regret Minimization (RRM) theory. Notably, RUM-based DCMs have been the go-to choice in the RM landscape. However, numerous pieces of research show that RRM-based DCMs are superior to RUM-based DCMs in scenarios in which customers tend to minimize regret.

This study contributes by formulating novel MILP models tailored for both capacitated and uncapacitated supply scenarios, with a specific focus on customers as regret minimizers based on RRM theory.\\
Particularly, the anticipated regret of all alternatives was computed through the pairwise comparison of all alternatives using MILP model constraints. Additionally, price-related regret included price level variables to compute regret for different price levels. Simultaneously, price level variables were integrated into the objective function of the MILP models to maximize revenue. By doing so, demand and supply were linked to both maximize revenue while accounting for customer regret minimization behavior. \\
The results not only validate the viability of these models but also underscore the substantial impact of varied customer behavior interpretations on revenue, pricing, and customer choices.

uture research endeavors may involve addressing large-scale datasets to solve RUM-based and RRM-based MILP models, shedding light on their real-world applicability. Specifically, the structures of MILP models indicate a significant number of variables and constraints, even for small instances.
Given the potential computational challenges, particularly with the heightened complexity of RRM-based MILP models, adopting Benders-type heuristics, column generation, or meta-heuristics emerges as a pragmatic solution.

Furthermore, acknowledging that both RUM-based and RRM-based MILPs often result in extreme values for prices, utility, and regret, the development of bi-objective MILPs is proposed. Introducing a second objective aimed at maximizing utility or minimizing regret could provide deeper insights. Exploring Pareto optimal frontiers, i.e., investigating the trade-offs between revenue or prices and utility or regret, in such models promises nuanced information about customer satisfaction, empowering suppliers with a more informed approach to decision-making in Revenue Management.

\section{Declaration of competing interest}

The author affirms that there are no identifiable conflicting financial interests or personal associations that might have seemed to impact the work described in this article.

\section{Data and Code availability}

Data and code are available upon  request.

\section{Acknowledgments}

We acknowledge Grenoble Alpes university IDEX scholarship. Dr. Pierre Lemaire and Dr. Iragael Joly provided important support.


\bibliographystyle{elsarticle-num} 
\bibliography{references}
\end{document}